\documentclass[11pt]{article}
\usepackage[margin=1in]{geometry}
\usepackage{amsmath,amssymb,amsthm}
\usepackage{tikz}
\usepackage{hyperref}
\usepackage{setspace}
\newtheorem{theorem}{Theorem}[section]
\newtheorem{proposition}[theorem]{Proposition}
\newtheorem{lemma}[theorem]{Lemma}

\theoremstyle{definition}
\newtheorem{definition}[theorem]{Definition}
\newtheorem{example}[theorem]{Example}
\theoremstyle{remark}
\newtheorem{remark}[theorem]{Remark}

\newcommand{\Sig}{\Sigma}
\newcommand{\X}{\mathcal X}
\newcommand{\Om}{\Omega}
\newcommand{\ParrySimplex}{\Delta_{\mathrm{P}}}

\title{The Parry Simplex: Freezing Phase Transitions for Topological Markov Chains\\}
\author{Shamsa Ishaq\thanks{shamsa.ishaq@lcwu.edu.pk}}
\date{}

\begin{document}
\maketitle

\begin{abstract}
We prove that the pressure function for the Hofbauer potential, defined via the distance to a topological Markov chain (subshift of finite type) $\X$ inside the one-sided full shift, exhibits a freezing phase transition. When $\X$ is irreducible (in particular when mixing), the post-transition equilibrium measure is unique, equal to the Parry measure of $\X$. When $\X$ is reducible with several communicating classes tied for maximal entropy and no chaining among them, the freezing transition persists, but the post-transition equilibrium measure ranges over the full simplex spanned by the Parry measures of the maximal-entropy components.\\


\end{abstract}

\section{Introduction}

\section*{Introduction}

In dynamical systems, a measure of maximal entropy (MME) is an invariant probability measure that maximizes the (measure theoretic) entropy among all invariant probability measures of a given system.
Thermodynamic formalism gives a way to study MME using ideas from physics. A dynamical system $(X,f)$ is considered as a physical  system, a continuous function $\varphi$ on the system as a potential, and 
the parameter $t > 0$ as an inverse temperature. By \cite{Walters}, the pressure at temperature $t$ is
\[
P(t) = \sup_{\mu} \left( h_\mu + t \int_X \varphi \, d\mu \right),
\]
where the supremum is over invariant probability measures $\mu$ and $h_\mu$ is the entropy of $\mu$. A measure attained this supremum is 
called an equilibrium state (measure) at temperature $t$ \cite{Ruelle2004,RuellePF}.\\
When $\varphi=0$ and the system is a mixing subshift of finite type, the equilibrium state is simply the unique measure of maximal entropy \cite{Parry}. More generally, if  $\varphi$ is H\"older continuous, there is always a unique
equilibrium state, and the pressure function is real analytic \cite{Ruelle2004,Bowen2008}. In other words, the system always has one perfect state at each temperature $t$.

A \emph{phase transition} appears when the pressure function $P(t)$ is not real analytic at some point $t_0$, which happens
exactly when there is more than one equilibrium state at $t_0$. To get a phase transition, Hofbauer \cite{Hofbauer} proved that one does not need to change the dynamics at all; it is enough to reduce the regularity of the potential below H\"older continuity. The potential in Hofbauer's case is defined on the full shift, with value depending on how close a point is to a single fixed point, and he proved that some critical point $t_0$ exists such that beyond $t_0$, the only equilibrium state is the Dirac measure at that fixed point. This is called a
\emph{freezing} phase transition: after $t_0$, the pressure stops changing, and the system \emph{freezes} onto one special state.

Later, Bruin and Leplaideur \cite{BruinLeplaideur2013,BL} generalized Hofbauer's work by using \emph{substitution subshifts} instead of a fixed point, the \emph{Fibonacci} and \emph{Thue-Morse} subshifts, which model quasicrystals. These are richer invariant sets than a single point: they have low complexity (zero entropy), but infinitely many points, and are uniquly ergodic. They proved  that a freezing phase transition happens for such Hofbauer-type potentials, with the pressure equal to the entropy of the
subshifts (which is zero).  However, the original proofs had a gap: a counting argument in the proof was not fully correct. Fixing this gap took a long time and a lot of work, and was only recently completed for the Thue-Morse case, using a refined combinatorial framework
that carefully tracks \emph{accident} times, moments where an orbit unexpectedly falls out of the subshift's language \cite{BCHL2025}. Even finding the exact critical temperature $t_0$ has proved difficult: only a partial estimate is known so far for Fibonacci case \cite{IshaqLeplaideur}.

This Hofbauer-type potential gives a natural motivation for defining a potential based on the distance to a general subshift of finite type  (mixing, irreducible, or reducible). Several questions then naturally arise: does the pressure function exhibit the same behavior, that is, 
does a freezing phase transition exist? Is the pressure function eventually equal to the topological entropy? And are the post-transition (frozen) states unique, coinciding with the unique measure of maximal entropy? Since subshifts of finite type are very well understood
their entropy, irreducible components, and periods can all be read  directly from a transition matrix using Perron Frobenius theory \cite{HJ,LM}, these questions are natural to pursue with the same tools used in the Bruin and Leplaideur case.

In this article we answer all three questions. Using the same inducing-scheme techniques as in \cite{Lep}, we show that the first two questions have a positive answer: a freezing phase transition always occurs, and the pressure does eventually equal the topological 
entropy. However, we show that the post transition measure is \emph{not} always unique(In case of reducible subshift). A similar non-uniqueness was found by Kucherenko and Thompson \cite{KucherenkoThompson2020} for suspension flows over the full shift, where a roof function built from a subshift of finite type with several tied maximal-entropy components produces several measures of maximal entropy for the flow, one per component. In our case, however, 
this non-uniqueness arises differently: it appears only as $t \to \infty$, as the outcome of a	 freezing transition, whereas in  \cite{KucherenkoThompson2020} it is not connected to temperature at all, it is built directly into the construction from the start.

To prove our result, we need to extend the classical estimate $T^n_{\beta\alpha}\sim C\eta^n$ from Perron-Frobenius theory to the reducible case, where paths between two states may pass through several maximal entropy components; getting this estimate right requires being careful about how we count the number of such components a path can pass through, which we explain in detail in Section~5. 

The rest of the article is organized as follows. Section~2 gives the background on words, symbolic dynamics, and thermodynamic formalism that we will need. Section~3 introduces the Hofbauer potential
built from a subshift of finite type, and states the main theorem. Section~4 proves the Theorem for mixing case. Section~5 develops the Perron-Frobenius tools needed for the reducible case, states and
proves the general theorem, and works through a numerical example showing exactly when this non-uniqueness does, and does not, happen.

\section{Basic Definitions and Examples}

Let $\mathcal S=\{0,1,\dots,m-1\}$ be a finite alphabet. A \emph{finite word}
is a finite sequence of symbols; if $u=u_0u_1\cdots u_{p-1}$, then $p=|u|$ is
its length, and $\epsilon$ denotes the empty word. Concatenation of words $u$
and $v$ is written $uv$. A finite word $u$ is a \emph{factor} of $v$ if
$v=wuw'$ for some words $w,w'$; it is a \emph{prefix} if $w=\epsilon$, a
\emph{suffix} if $w'=\epsilon$, and an \emph{inner factor} if $w,w'\ne
\epsilon$. We denote by $\Sig^+$ the set of finite words over $\mathcal S$.

A (one-sided) \emph{infinite word} is $z=z_0z_1z_2\cdots$ with $z_i\in
\mathcal S$; $\Sig=\mathcal S^{\mathbb N}$ is the set of all such words. The
notions of prefix, factor, and suffix extend naturally to infinite words.
Given $z\in\Sig$, its \emph{language} $L(z)$ is the set of all its finite
factors, and $L_n(z)=\{u\in L(z):|u|=n\}$.

Given $x,y\in\Sig$, their distance is $d(x,y)=2^{-\ell(x,y)}$, where
$\ell(x,y)$ is the length of the longest common prefix ($d(x,y)=0$ if
$x=y$). The \emph{shift} $\sigma:\Sig\to\Sig$ is defined by
$\sigma(x_0x_1x_2\cdots)=x_1x_2x_3\cdots$; it is Lipschitz continuous. $(\Sig,\sigma)$ is the \emph{one-sided full shift} over
$\mathcal S$.

For a finite word $w=w_0\cdots w_{n-1}$, the \emph{cylinder} $[w]=\{y\in
\Sig: y_i=w_i,\ 0\le i\le n-1\}$. A finite word $u$ is a \emph{return word}
of the cylinder $J=[w]$ if: (1) $w$ is a prefix of $uw$; (2) $w$ is not an
inner factor of $uw$; (3) $\min\{k:\sigma^k(ux)\in[w]\}=|u|$ for some
$x\in[w]$. We write $R_J$ for the set of return words to $J$; if $y=ux\in J$
for $u\in R_J$, then $|u|$ is the \emph{first return time}.

A subset of $\Sig$ that is closed and $\sigma$-invariant is a
\emph{subshift}. A subshift $\mathcal{X}$ is a \emph{subshift of finite type}
(SFT), or topological Markov chain (TMC), if there is a $0/1$ transition matrix
$T=(t_{ij})$ of order $m\times m$ such that $x\in\mathcal{X}$ iff
$t_{x_ix_{i+1}}=1$ for all $i$; we write $\mathcal{X_T}$ for this subshift or $\mathcal{X}$ is subshift with transition matrix $T$. The language of subshift $\mathcal{X}$, is denoted by
$L(\mathcal{X})$ is the set of finite words appearing in elements of
$\mathcal{X}$, and $L_n(\mathcal{X})=\{u\in L(\mathcal{X}):|u|=n\}$.

An SFT $\mathcal{X}$ with associated matrix $T$ is \emph{irreducible} if for every $u,v\in L(\mathcal{X})$ there
is $w\in L(\mathcal{X})$ with $uwv\in L(\mathcal{X})$; it is \emph{mixing} if for every
$u,v$ there is $N_{u,v}$ such that for every $n\ge N_{u,v}$ some $w\in
L_n(\mathcal{X})$ satisfies $uwv\in L(\mathcal{X})$. Equivalently, in term of the transition matrix T: 
$\mathcal{X}$ is irreducible if and only if T is irreducible matrix (for every $i,j$ there is an $n\geq 1$ with $T^{n}_{ij}>0)$, 
and $\mathcal{X}$ is mixing if and only if $T$ is primitive matrix (there is a single $n$ with $T^{n}_{ij}>0)$, for all $i, j$ at once). A SFT that is not irreducible is called reducible, corresponding to a reducible transition matrix.

The entropy of $\mathcal{X}$ is $h_\sigma=\lim_{n\to\infty}\frac1n \log
|L_n(\mathcal{X})|$ (the limit always exists). Given an adjacency matrix $T$ of
size $m$, $|L_n(\mathcal{X})|=\sum_{i,j=1}^m T^n_{ij}$, and $h_\sigma=\log\eta$
where $\eta$ is the dominating eigenvalue of $T$ \cite[Ch.~4]{LM}. For a
mixing SFT there is $n\ge1$ with $T^n>0$ \cite{Bowen2008}.

\begin{theorem}[Perron]
\label{thm:perron}
Let $T$ be a positive matrix of order $n$. Then: (i) $\rho>0$ is an
algebraically simple eigenvalue of $T$; (ii) there is a unique positive
vector $u$, $\|u\|=1$, with $Tu=\rho u$; (iii) there is a unique positive
vector $v$, $\|v\|=1$, with $v^TT=\rho v^T$ and $\langle u,v\rangle=1$;
(iv) $(\rho^{-1}T)^m\to uv^T$ as $m\to\infty$.
\end{theorem}

Let $M(\Sig)$ be the Borel probability measures on $\Sig$, $C(\Sig)$ the
continuous complex-valued functions on $\Sig$. A measure $\mu$ is
$\sigma$-invariant if $\mu(\sigma^{-1}(B))=\mu(B)$ for all Borel measurable $B$; let
$M_\sigma(\Sig)$ be the set of such measures, a nonempty compact convex
subset of $M(\Sig)$ \cite{Bowen2008}. An SFT is \emph{uniquely ergodic} if it
carries exactly one invariant probability measure.

For $n\in\mathbb N$, let $C_n=\{[w]:w\in\Sig^+_n\}$; there are $|\mathcal
S|^n$ such cylinders. The \emph{Kolmogorov entropy} of $\mu$ is
\begin{equation}
h_\mu=\lim_{n\to\infty}\frac{-1}{n}\sum_{[w]\in C_n}\mu[w]\log\mu[w],
\end{equation}
the limit always existing \cite[Lemma~1.19]{Bowen2008}. 

For an SFT $\mathcal{X}$, a measure $\mu \in M_\sigma(\mathcal{X})$ with $h_\mu = h_\sigma$ is called a
\emph{measure of maximal entropy}; we write
\[
\mathrm{MME}(\mathcal{X}) \;:=\; \{\mu \in M_\sigma(\mathcal{X}) : h_\mu = h_\sigma\}
\]
for the set of all measures. When $\mathcal{X}$ is irreducible, intrinsic ergodicity \cite{Parry}
guarantees that $\mathrm{MME}(\mathcal{X})$ is a singleton, and its unique element is called the
\emph{Parry measure} of $\mathcal{X}$, denoted $\mu_{\mathcal{X}}$; in this case the two notions coincide.
When $\mathcal{X}$ is reducible with maximal entropy components $\mathcal{X}_1,\dots,\mathcal{X}_k$, the Parry measures
$\mu_{_{\mathcal{X}_1}},\dots,\mu_{_{\mathcal{X}_k}}$ are the extreme points of the larger simplex;
\[
\mathrm{MME}(\mathcal{X}) \;=\; \Big\{\, \textstyle\sum_{i=1}^k c_i\,\mu_{\mathcal{X_i}} \;:\; c_i \ge 0,\ \sum_i c_i = 1 \,\Big\},
\]
we call the above simplex, the \emph{Parry simplex} of $\mathcal{X}$ and denoted it by;

\[
\mathrm{MME}(\mathcal{X}) = \ParrySimplex{\mathcal{X}} ,
\]
i.e, the simplex spanned by the Parry measures of the maximal components.  which is exactly the set of post transition equilibrium measures identified in Section 5.

A continuous $\varphi:\Sig\to\mathbb R$ is a \emph{potential function}; its
\emph{pressure} over $(\Sig,\sigma)$ is
\[
P(\varphi)=\sup_{\mu\in M_\sigma(\Sig)}\Big(h_\mu+\int\varphi\,d\mu\Big).
\]
A measure $\mu^+$ realizing the supremum is an \emph{equilibrium measure}
\cite{Walters}. Since $\mu\mapsto h_\mu$ is upper semicontinuous on the
compact space $M_\sigma(\Sig)$, every continuous $\varphi$ admits at least
one equilibrium measure \cite{Bowen2008}. The \emph{pressure function}
$t\mapsto P(t):=P(t\varphi)$ is convex; where $P'(t)$ exists and $\mu_t$ is
the equilibrium measure for $t\varphi$, $P'(t)=\int\varphi\,d\mu_t$, and as
$t\to\infty$ the graph of $P$ has an asymptote of slope $m=\max\{\int
\varphi\,d\mu:\mu\in M_\sigma(\Sig)\}$ \cite{BLL}.

Ruelle's Perron--Frobenius theorem \cite{RuellePF} gives, for H\"older
$\varphi$, a unique equilibrium measure and real-analytic
pressure function. However, when multiple equilibrium measures exists at some $t_0$, then the pressure function $P(t)$ may
lose analyticity: assume that $\mu^+,\mu^-$ are two equilibrium states for $t_0\varphi$ and
$\int\varphi\,d\mu^+\ne\int\varphi\,d\mu^-$, then $P(t)$ is not 
differentiable at $t_0$. We say $P(t)$ has a \emph{phase transition} at $t_0$
if $P(t)$ is not real analytic at $t_0$, and the transition is \emph{freezing} at $t_0$, if
$P(t)$ becomes affine after $t_0$.

\section{Setting, Main Theorem, and Idea of Proof}

Let $\mathcal{X}$ be a subshift of finite type in the one-sided full shift $\Sigma$, with language
$L(\mathcal{X})$. For $x = x_0x_1\cdots \in \Sigma$, set
\[
\delta(x) := \max\{n : \forall\, k\le n,\ x_0\cdots x_k \in L(\mathcal{X})\} \;\le\; +\infty;
\]
thus $\delta(x)=+\infty$ iff $x\in\mathcal{X}$. Similarly for $u\in\Sigma^+\setminus L(\mathcal{X})$,
$\delta(u)$ is the length of the longest common prefix of $u$ in $L(\mathcal{X})$. By definition,
$d(x,\mathcal{X}) = 2^{-\delta(x)-1}$.

Fix $w\in\Sigma^+\setminus L(\mathcal{X})$, $J=[w]$; then $\delta(x)=\delta(w)$ for all $x\in J$.
Fix $N\gg\delta(w)$ and $A>0$; define the \emph{Hofbauer potential}
\[
\varphi(x) = \begin{cases} -\log\big(1+\tfrac{1}{\delta(x)}\big) & \delta(x)\ge N,\\[2pt] -A & \delta(x)<N, \end{cases} \tag{2}
\]
so that $\varphi\equiv0$ on $\mathcal{X}$.

\begin{theorem}[Main Theorem]\label{thm:main}
Let $\mathcal{X}$ be a subshift of finite type with entropy $\xi=\log\eta$ in $\Sigma$, and
$\varphi$ the Hofbauer potential above. Assume $\mathcal X$ has non-chained maximal-entropy
components $\mathcal{X}_1,\dots,\mathcal{X}_k$ ($k\ge1$; see Section~5.1). There is a transition point $t_0>0$ such that:
\begin{enumerate}
\item For $0<t<t_0$ there is a unique equilibrium measure $\mu_t$, of full support; $P(t)$ is real
analytic and $P(t)>\xi$.
\item For $t>t_0$, $P(t)=\xi$, and
\[
\ParrySimplex{\mathcal{X}} =\; \Big\{\, \textstyle\sum_{i=1}^k c_i\,\mu_{_{\mathcal{X}_i}} \;:\; c_i\ge0,\
\sum_i c_i=1 \,\Big\}
\]
is exactly the set of equilibrium measures for $t\varphi$: The transition is freezing, and the post-transition equilibrium measure
is unique equal to the Parry measure $\mu_{\mathcal X}$  if and only if $k=1$, which holds in particular whenever $\mathcal X$ is irreducible or mixing.
\end{enumerate}
\end{theorem}

\begin{remark}
When $\mathcal X$ is \emph{mixing} (the leading case, $k=1$ trivially), Theorem~\ref{thm:main} is
proved directly in Section~4, and
$\mathrm{MME}(\mathcal X)=\{\mu_{\mathcal X}\}$. Section~5 extends the argument to $\mathcal X$
irreducible but not mixing (Proposition~5.7, still $k=1$) and then to $\mathcal X$ reducible with
$k\ge2$ non-chained maximal components (Theorem~5.10), where $\mathrm{MME}(\mathcal X)$ becomes a
 simplex rather than a single point, the non-uniqueness phenomenon after post-transition, that is the main
contribution of this article.
\end{remark}

To prove Theorem~\ref{thm:main} in the mixing case we follow Leplaideur's inducing-scheme method
\cite{Lep}. For $t\ge0$, $z\in\mathbb{R}$, $x\in J$, set
\[
\lambda_{t,z} = \sum_{u\in R_J} e^{tS_{|u|}\varphi(ux) - |u|z}, \tag{3}
\]
where $S_{|u|}\varphi(y) = \sum_{i=0}^{|u|-1}\varphi(\sigma^i(y))$ is the Birkhoff sum along the
orbit segment $O_u(y) = \{\sigma^i(y) : 0\le i\le |u|-1\}$. The general (reducible) case is treated
in Section~5 by replacing the scalar Perron-Frobenius asymptotic (Lemma~5.4).

From \cite{Lep},  For each fixed $t\ge0$ there is a minimal critical number $z_c(t)\ge-\infty$
such that for $z>z_c(t)$, $\lambda_{t,z}<+\infty$ for all $x\in J$, and
$z_c(t)\le P(t)$. Moreover
\[
z_c(t)=P(\Sig_J,t)=\sup_{\mu\in M_\sigma(\Sig_J)}\Big(h_\mu+t\int\varphi\,
d\mu\Big),
\]
where $\Sig_J=\{x\in\Sig:\sigma^n(x)\notin J\ \forall n\in\mathbb N\}$, i.e.\
$z_c(t)$ is the pressure of $t\varphi$ restricted to points whose orbit never
meets $J$.

\section{Proof of Theorem~\ref{thm:main} (The Mixing case)}
\label{sec:mixing-proof}

Let $\X$ be a mixing SFT with transition matrix $T$; there is $N_\X$ with
$T^n>0$ for $n\ge N_\X$. Fix $N\gg N_\X$ and $\varphi$ as in~(2). For
$\alpha,\beta\in\mathcal S$ with $\alpha\beta\notin L(\X)$, set $J=[\alpha
\beta]$; then $\varphi\equiv-A$ on $J$. Let $\mathcal F_\X=\{x\in\Sig_{
\mathcal S}:\delta(x)\le N-1\}$ (the \emph{free region}) and $\mathcal E_\X=
\{x\in\Sig_{\mathcal S}:\delta(x)\ge N\}$ (the \emph{excursion region}); clearly,
$J\subset\mathcal F_\X$ and $\X\subset\mathcal E_\X$.

For $u\in R_J$, $i\in[[1,|u|-1]]$ is an \emph{accident time} if
$\delta(\sigma^i(ux))>\delta(\sigma^{i-1}(ux))-1$; otherwise
$\delta(\sigma^i(ux))=\delta(\sigma^{i-1}(ux))-1$.

\begin{figure}[h]
\centering
\begin{tikzpicture}[scale=1.0]
\draw[->] (0,0) -- (10,0) node[right] {};
\draw[thick] (0.5,0) -- (0.5,0.25) node[above] {\footnotesize $\alpha\beta$};
\draw[thick] (9.5,0) -- (9.5,0.25) node[above] {\footnotesize $\alpha\beta$};
\draw[thick] (3.0,0) -- (3.0,0.25) node[above] {\footnotesize $\alpha^{(1)}\beta^{(1)}$};
\draw[thick] (6.5,0) -- (6.5,0.25) node[above] {\footnotesize $\alpha^{(k-1)}\beta^{(k-1)}$};
\draw[<->] (0.6,-0.4) -- (2.9,-0.4) node[midway,below] {\footnotesize $w^{(1)}$};
\draw[<->] (3.1,-0.4) -- (6.4,-0.4) node[midway,below] {\footnotesize $\cdots$};
\draw[<->] (6.6,-0.4) -- (9.4,-0.4) node[midway,below] {\footnotesize $w^{(k)}$};
\draw[<->,dashed] (0.5,-1.1) -- (3.0,-1.1) node[midway,below]
  {\footnotesize safe zone};
\draw[<->,dashed] (3.0,-1.1) -- (6.5,-1.1) node[midway,below]
  {\footnotesize safe zone};
\draw[<->,dashed] (6.5,-1.1) -- (9.5,-1.1) node[midway,below]
  {\footnotesize safe zone};
\end{tikzpicture}
\caption{Decomposition~(4) of a return word with $k$ accidents.}
\label{fig:return-decomp}
\end{figure}

\begin{lemma}
\label{lem:accident}
Let $u$ be a return word to $J=[\alpha\beta]$ and $k\in[[1,|u|-1]]$ an
accident time with $\delta(\sigma^k(ux))=m$ for $x\in J$. Then for all
$i\in[[k,k+m]]$, $\delta(\sigma^i(ux))=\delta(\sigma^{i-1}(ux))-1$.
\end{lemma}

\begin{proof}
Let $u_k\cdots u_{k+m}\in L(\X)$, $u_k\cdots u_{k+m}u_{k+m+1}\notin L(\X)$, so
$u_{k+m}u_{k+m+1}=\alpha^1\beta^1$ is a forbidden word. If some
$i\in[[k,k+m]]$ had $\delta(\sigma^i(ux))>\delta(\sigma^{i-1}(ux))-1$, then
$\delta(\sigma^i(ux))>m-i$, forcing $u_{k+m}u_{k+m+1}$ to appear inside a
word of $L(\X)$ of length $\delta(\sigma^i(ux))$ -- a contradiction.
\end{proof}

Consequently a return word with $k$ accidents decomposes as
\[
u=\alpha\beta\, w^{(1)}\alpha^{(1)}\beta^{(1)} w^{(2)}\alpha^{(2)}\beta^{(2)}
\cdots \alpha^{(k-1)}\beta^{(k-1)} w^{(k)} \alpha\beta, \tag{4}
\]
where $\beta^{(s-1)}w^{(s)}\alpha^{(s)}\in L(\X)$ and $\alpha^{(s)}
\beta^{(s)}\notin L(\X)$ for $s=1,\dots,k$ (with $\alpha^{(k)}=\alpha$,
$\beta^{(0)}=\beta$); see Figure~\ref{fig:return-decomp}.

Let $\mathcal O^+(u)=\{\sigma^i(ux):0\le i\le|u|-1,\ x\in J\}$. A return word
$u$ is of \emph{type~1} if $\mathcal O^+(u)\cap\mathcal E_\X=\varnothing$; of
\emph{type~2} if the orbit enters $\mathcal E_\X$ exactly once; of
\emph{type~3} if it enters multiple times. Write $\mathcal T_1,\mathcal
T_2,\mathcal T_3$ for the corresponding sets of return words, and $\Om$ for
the set of forbidden blocks of $\X$.

\subsection*{Case $|\Om|=1$}

The return word has a single accident, so $R_J=\{\alpha\beta w: w\in L(\X),
\ \beta w\alpha\in L(\X)\}$. Since $\X$ is mixing, for each $n\ge N$ there is
$w\in L_n(\X)$ with $\alpha w\beta\in L(\X)$ giving a type-2 return word; if
$|u|<N$ the return word is type~1. There is no type-3 return word.

Splitting $\lambda_{t,z}=\lambda_{t,z}(\mathcal T_1)+\lambda_{t,z}(\mathcal
T_2)$ at $n=N$: for $u\in R_J$, $|u|\le N-1$, $S_{|u|}\varphi(ux)-|u|z=
-(A+z)|u|$, so
\[
\lambda_{t,z}(\mathcal T_1)=\sum_{n=1}^{N-1} c_n e^{-t(A+z)n},
\]
$c_n$ the multiplicity of type-1 words of length $n$. For $|u|\ge N$,
$S_{|u|}\varphi(ux)=-NA+\log N-\log n$ for $x\in J$, so
\[
\lambda_{t,z}(\mathcal T_2)=e^{-t(NA-\log N)}\sum_{n\ge N} d_n\Big(
\frac1{n^t}\Big)e^{-nz},
\]
$d_n$ the multiplicity of type-2 words of length $n\ge N$. Since $N\gg
N_\X$, the number of length-$n$ paths from $\beta$ to $\alpha$ is
$T^n_{\beta\alpha}$, and Theorem~\ref{thm:perron} gives $C\eta^n\approx d_n$
($n\ge N$), $\eta$ the spectral radius of $T$, $C$ depending on
$(\beta,\alpha)$. Hence
\[
\lambda_{t,z}(\mathcal T_2)=Ce^{-t(NA-\log N)}\sum_{n\ge N}\Big(\frac\eta
{e^z}\Big)^n \frac1{n^t}. \tag{5}
\]
This series converges iff $\eta/e^z<1$, i.e.\ $z\ge\xi:=\log\eta$, the
entropy of $\X$. At $z=\xi$, 
\[
\lambda_{t,\xi}=\sum_{n=1}^{N-1}c_n e^{-tAn}+e^{-tC(A,N)}\sum_{n\ge N}
\frac1{n^t}, \tag{6}
\]
converging for all $t\ge0$, so $z_c(t)=\xi$, with $c_n\ge0$ and
$C(A,N)=NA-\log N>0$.

Since $\frac{d}{dt}\lambda_{t,\xi}\le0$, $\lambda_{t,\xi}$ is decreasing;
$\lambda_{t,\xi}\to+\infty$ as $t\to1$ (via $\lambda_{t,z}(\mathcal T_2)$),
and $\lambda_{t,\xi}\to0$ as $t\to\infty$. By continuity there is $t_c>1$
with $\lambda_{t_c,\xi}=1$. By \cite{Lep}[Theorem 4]: for $t<t_c$,
$\lambda_{t,z}>1$ for all $z\ge\xi$, and since $z\mapsto\lambda_{t,z}$ is
decreasing, there is $z(t)>\xi$ with $P(t)=z(t)$, with a unique equilibrium
measure $\mu_t$ and $P$ real analytic; for $t>t_c$, $\lambda_{t,z}<1$ for
all $z\ge\xi$, so no equilibrium measure gives positive weight to $J$, and
$P(t)=\xi$.

\subsection*{Case $|\Om|\ge2$}

Now a return word with $k$ accidents has the general form~(4). Let
$\mathcal S_J\subset R_J$ be the return words whose orbit enters
$\mathcal E_\X$ exactly once, and $\lambda^{S_J}_{t,z}=\sum_{u\in\mathcal
S_J}\sum_{|u|=n} e^{tS_{|u|}\varphi(ux)-|u|z}$. Writing $u=FEF$ (free part
$F$, excursion part $E$),
\[
\lambda^{S_J}_{t,z}\le \big(F(t,z)\big)^2 E(t,z),
\]
where $F(t,z)=\sum_F\sum_{|F|=n} e^{tS_{|F|}\varphi(Fx)-|F|z}\le
\sum_{n\ge0}(me^{-tA-z})^n$ converges for $z>\log m-tA$ ($m=|\mathcal S|$),
and (by the same computation as~(5))
\[
E(t,z)=Ce^{-t(NA-\log N)}\sum_{n\ge N}\Big(\frac\eta{e^z}\Big)^n\frac1{n^t}
\]
converges for $z>\log\eta$. At $z=\xi$,
\[
F(t,\xi)\le \sum_{n\ge1}\Big(\frac m\eta e^{-tA}\Big)^n=\frac m{\eta e^{tA}-m},
\qquad E(t,\xi)=Ce^{-t(NA-\log N)}\zeta(t).
\]


Summing over return words entering $E_{\mathcal X}$ exactly $n$ times gives a contribution
$\le F(t,\xi)^{n+1}E(t,\xi)^n$. Since $F(t,\xi)$ and $E(t,\xi)$ are already known explicitly, we sum
this geometric series directly rather than passing through a further bound: whenever
$F(t,\xi)E(t,\xi)<1$,
\[
\lambda_{t,\xi} \;\le\; \sum_{n\ge0} F(t,\xi)^{n+1}E(t,\xi)^{n}
\;=\; \frac{F(t,\xi)}{1-F(t,\xi)E(t,\xi)}. \tag{7}
\]
Substituting the explicit expressions for $F(t,\xi)$ and $E(t,\xi)$, and writing
$\phi := \dfrac{m}{\eta e^{tA}-m}$ for brevity,
\[
\lambda_{t,\xi} \;\le\; \frac{\phi}{1-\phi\, Ce^{-t(NA-\log N)}\zeta(t)},
\]
valid precisely when $\phi\, Ce^{-t(NA-\log N)}\zeta(t) < 1$. Together with the lower bound
\[
\lambda_{t,\xi} \;\ge\; C\, e^{-t(NA-\log N)}\,\zeta(t),
\]
obtained exactly as before from the type-2 (single-excursion) return words alone.


As $t\to1$, $\zeta(t)\to+\infty$ so $\lambda_{t,\xi}\to+\infty$ by the lower bound (since $\phi$
stays bounded away from $0$); as $t\to\infty$, both $F(t,\xi),E(t,\xi)\to0$, so
$\lambda_{t,\xi}\to0$. By continuity and monotonicity there is $t_c>1$ with $\lambda_{t_c,\xi}=1$,
and the conclusion of Theorem~3.1 follows as in the case $|\Omega|=1$. This completes the proof of
Theorem~3.1. $\blacksquare$

\section{Phase Transition for General Subshifts of Finite Type}
\label{sec:general}

We now remove the mixing hypothesis on $\X$ entirely. The combinatorial
framework of Section~\ref{sec:mixing-proof}  Lemma~\ref{lem:accident}, the
decomposition~(6), and the classification of return words into
$\mathcal T_1,\mathcal T_2,\mathcal T_3$  depends only on the language
$L(\X)$ and the function $\delta$, never on mixing or irreducibility, and
carries over unchanged. Only the scalar asymptotic $T^n_{\beta\alpha}\sim
C\eta^n$, which used Theorem~\ref{thm:perron}(iv), needs to be replaced.

\subsection{Reducible Perron-Frobenius normal form}

Let $T$ be the $m\times m$ transition matrix of $\X$ on alphabet
$\mathcal S$. After permuting states, $T$ is in block upper-triangular
\emph{normal form}
\[
T=\begin{pmatrix}
T_{11} & \ast & \cdots & \ast \\
0 & T_{22} & \cdots & \ast \\
\vdots & & \ddots & \vdots \\
0 & 0 & \cdots & T_{rr}
\end{pmatrix},
\]
where each diagonal block $T_{ii}$ is irreducible and corresponds to a
communicating (strongly connected) class $C_i$ of the underlying directed
graph, $i=1,\dots,r$ \cite[Ch.~8.3]{HJ}. Write $\eta_i=\rho(T_{ii})$ and
$\xi_i=\log\eta_i$, the spectral radius and entropy of each class $C_i$. Also 
\[
\eta=\max_{1\le i\le r}\eta_i,\qquad \xi=\log\eta=h_\sigma(\X),
\]
the entropy of $\X$ \cite[Ch.~4]{LM}. Each irreducible class $C_i$ has a
period $p_i\ge1$ ($p_i=1$ iff $\sigma|_{\Sig_{C_i}}$ is mixing) and a unique
Parry measure $\mu_{C_i}$ by intrinsic ergodicity of irreducible SFTs
\cite{Parry}.

Let $\mathcal G$ be the \emph{condensation}: the DAG on $\{1,\dots,r\}$ with
an edge $i\to j$ whenever some state of $C_i$ has an edge into $C_j$ in the
original graph, $C_i\ne C_j$. Write $C_i\succeq C_j$ if there is a directed
path (possibly of length $0$) from $i$ to $j$ in $\mathcal G$; this is a
partial order since $\mathcal G$ is acyclic. Set
\[
I_{\max}=\{i:\eta_i=\eta\},\qquad \X_{\max}=\bigcup_{i\in I_{\max}}\X_i,
\]
where $\X_i\subset\X$ is the irreducible sub-SFT generated by $C_i$.

\begin{definition}[Chain length]\label{def:5.1}
For $\beta\in C_i$, $\alpha\in C_j$ with $C_i\succeq C_j$, let $\mathcal P(i,j)$ denote the set of
directed paths $P=(s_1,\dots,s_\ell)$ in the condensation DAG $\mathcal{G}$ with $s_1=i$, $s_\ell=j$. Since
$\mathcal{G}$ is a finite acyclic graph, $\mathcal P(i,j)$ is a \emph{finite}, nonempty set. Define
\[
d(\beta,\alpha):=\max_{P\in\mathcal P(i,j)}\ \#\{a: C_{s_a}\in I_{\max}\},
\]
the largest number of maximal classes encountered along any \emph{single} directed route from
$C_i$ to $C_j$ in $\mathcal{G}$. If $C_i\not\succeq C_j$, $T^n_{\beta\alpha}=0$ for all $n$ and we set
$d(\beta,\alpha)$ undefined (irrelevant).
\end{definition}

\begin{definition}[Non-chained maximal classes]\label{def:5.2}
$\mathcal X$ has \emph{non-chained maximal classes} if the restriction of $\succeq$ to $I_{\max}$
is trivial: no two distinct $i,j\in I_{\max}$ satisfy $C_i\succeq C_j$.
\end{definition}

\begin{remark}\label{claim:dle1}
If $\mathcal X$ has non-chained maximal classes, then $d(\beta,\alpha)\le 1$ for \emph{every}
$\beta\in C_i$, $\alpha\in C_j$ with $C_i\succeq C_j$.

\end{remark}


\begin{lemma}[Reducible Perron--Frobenius asymptotics]\label{lem:5.3}
Let $\beta\in C_i$, $\alpha\in C_j$, $C_i\succeq C_j$, and let $d=d(\beta,\alpha)$ be as in
Definition~\ref{def:5.1}. If $d\ge 1$ there is $C=C(\beta,\alpha)>0$ with
\[
T^n_{\beta\alpha}\ \sim\ C\,n^{d-1}\eta^n \qquad (n\to\infty).
\]
If every $P\in\mathcal P(i,j)$ contains no element of $I_{\max}$ (i.e.\ $d=0$), then
$T^n_{\beta\alpha}=O(\rho^n)=o(\eta^n)$, where $\rho<\eta$ is the largest spectral radius occurring
among the classes visited by paths in $\mathcal P(i,j)$.
\end{lemma}

\begin{proof}
Any length-$n$ walk(return word of length $n$) from $\beta$ to $\alpha$ in the original graph visits a sequence of classes;
collapsing consecutive repeats yields a well-defined directed path in $\mathcal{G}$ from $i$ to $j$, called
the walk's \emph{shape}. Since $\mathcal P(i,j)$ is finite (Definition~\ref{def:5.1}), every walk
has a shape $P\in\mathcal P(i,j)$, and
\[
T^n_{\beta\alpha}=\sum_{P\in\mathcal P(i,j)} W_n(P),
\]
where $W_n(P)$ counts length-$n$ walks of shape exactly $P$.

Fix $P=(s_1,\dots,s_\ell)\in\mathcal P(i,j)$ and set $d_P:=\#\{a:C_{s_a}\in I_{\max}\}$. Split a
length-$n$ walk of shape $P$ into segments lying in the
successive classes $C_i=C_{s_1}\succ\cdots\succ C_{s_\ell}=C_j$ of the
chain, of lengths $k_1,\dots,k_\ell$ with $k_1+\cdots+k_\ell=n-O(1)$, the $O(1)$ absorbing the finitely many fixed transition edges used
to hop between consecutive classes of $P$ (these do not scale with $n$).


Apply the ordinary Perron--Frobenius theorem (Theorem~2.1) to each irreducible block $T_{s_as_a}$:
\begin{enumerate}
\item if $C_{s_a}\in I_{\max}$ (spectral radius $\eta$), a segment of length $k_a$ inside it
contributes $\Theta(\eta^{k_a})$;
\item if $C_{s_a}\notin I_{\max}$ (spectral radius strictly less than $\eta$), a segment of length
$k_a$ contributes only $o(\eta^{k_a})$, and is absorbed into a bounded multiplicative constant.
\end{enumerate}

Only the $d_P$ segments lying in maximal classes contribute genuine $\eta^{k_a}$-growth in $P$;
multiplying across segments and summing over all valid splittings gives, up to absorbed constants,
\[
W_n(P)=\Theta(\eta^n)\cdot\#\Big\{(k_1,\dots,k_{d_P}): \textstyle\sum k_a=n-O(1),\ k_a\ge 0\Big\}.
\]
This uses that the product $\prod_a\Theta(\eta^{k_a})=\Theta(\eta^{\sum k_a})=\Theta(\eta^{\,n-O(1)})=\Theta(\eta^n)$
depends only on the total $\sum k_a$, not on the individual splitting.

The number of ways to write $n+O(1)$ as an ordered sum of $d_P$ nonnegative integers is
\[
\binom{n+O(1)}{d_P-1}\ \sim\ C_{P}\,n^{d_P-1}\qquad(n\to\infty),\qquad C_{P}=\frac{1}{(d_P-1)!},
\]
by the standard stars-and-bars identity, valid for $d_P\ge 1$; if $d_P=0$ there is no maximal
segment and $W_n(P)=O(\rho_{_P^n})$ with $\rho_{_P}<\eta$ the largest spectral radius occurring on $P$.
Combining with $W_n(P)$ we get,
\[
W_n(P)\sim C_P\,n^{d_P-1}\eta^n\quad(d_P\ge1),\qquad W_n(P)=O(\rho_P^n)=o(\eta^n)\quad(d_P=0).
\]

Since $\mathcal P(i,j)$ is finite, summing preserves the order of the dominant term:
\[
T^n_{\beta\alpha}=\sum_{P\in\mathcal P(i,j)}W_n(P)\ \sim\
\Big(\sum_{P:\,d_P=d}C_P\Big)\,n^{d-1}\eta^n,\qquad d:=\max_{P\in\mathcal P(i,j)}d_P=d(\beta,\alpha),
\]
because every $P$ with $d_P<d$ contributes a strictly smaller power of $n$ (or is $o(\eta^n)$
altogether), and a finite sum cannot let a lower-order term overtake the leading one. Setting
$C:=\sum_{P:d_P=d}C_P>0$ gives $T^n_{\beta\alpha}\sim C\,n^{d-1}\eta^n$, as claimed. If $d=0$ (no
shape contains a maximal class), the same finite sum gives $T^n_{\beta\alpha}=O(\rho^n)=o(\eta^n)$
with $\rho=\max_P \rho_P<\eta$.
\end{proof}

\begin{remark}[Consistency check]\label{rem:5.4}
When $\mathcal X$ is mixing there is only one class ($r=1$), it is the unique maximal
class ($I_{\max}=\{1\}$), $\mathcal P(1,1)$ consists of the single trivial path, and $d\equiv1$.
Lemma~\ref{lem:5.3} then gives $T^n_{\beta\alpha}\sim Cn^0\eta^n=C\eta^n$, recovering exactly the
classical Perron--Frobenius estimate used in the mixing-case proof of Theorem~3.1. So
Lemma~\ref{lem:5.3} is a strict generalization of mixing case in section \ref{sec:mixing-proof}.
\end{remark}

\begin{example}\label{ex:5.5}
Take two disjoint golden-mean shifts on $\{0,1\}$ and $\{2,3\}$ (excluding ``$11$'', resp.\
``$33$''), entropy $\xi=\log\varphi$, $\varphi=(1+\sqrt5)/2$. Gluing them by a single one-way edge
$1\to2$ (no edge back) gives a condensation DAG that is a simple path with no branching, so
$\mathcal P(i,j)$ has exactly one element and $d(\beta,\alpha)=2$ unambiguously (both the original
formulation of Definition~\ref{def:5.1} and the corrected one agree here, since there is no
branching to disambiguate). A direct computation of $T^n$ confirms
\[
\frac{T^n_{0,2}}{n\varphi^n}\longrightarrow \frac15,\qquad \frac{T^n_{0,2}}{\varphi^n}\to\infty\quad(n\to\infty),
\]
the predicted linear correction. Chaining three such shifts $A\to B\to C$ gives
$T^n_{\beta\alpha}/(n^2\varphi^n)\to\text{const}$, confirming $d-1=2$. By contrast, feeding two
golden-mean shifts from a common transient state with \emph{no} edge between them ($d\equiv1$,
non-chained, Definition \ref{def:5.2}) gives $T^n_{\beta\alpha}/\varphi^n\to\text{const}$
with no polynomial correction.
\end{example}

\subsection{Proof of Theorem~\ref{thm:main} (The irreducible case)}


\begin{proposition}
\label{prop:A}
Theorem~\ref{thm:main} holds, with $\X$'s Parry measure as the
unique post-transition equilibrium measure, when $\X$ is any irreducible
subshift of finite type, mixing or not.
\end{proposition}

\begin{proof}
Only the estimate $d_n\approx C\eta^n$ preceding~(7) used mixing. If $\X$ is
irreducible with period $p$, its states split into cyclic classes
$S_0,\dots,S_{p-1}$, $T$ mapping $S_k$ into $S_{k+1\bmod p}$, and
$T^n_{\beta\alpha}>0$ only for $n$ in a single residue class $r(\beta,
\alpha)\bmod p$; along that residue class, primitivity of $T^p$ restricted
to a cyclic class gives $T^n_{\beta\alpha}\sim C\eta^n$ (the case $d=1$ of
Lemma~\ref{lem:5.3}. Restricting the sums in~ (5)--(7) to that
arithmetic progression, $\sum_{n\ge N,\,n\equiv r_0(p)}n^{-t}(\eta/e^z)^n$
converges under the same condition $z\ge\xi$ as the unrestricted sum, so
equations~(8) and the monotonicity argument for $\lambda_{t,\xi}$ carry
over unchanged. Uniqueness of the Parry measure for irreducible (not
necessarily mixing) SFTs is classical \cite{Parry}.
\end{proof}

\subsection{Proof of Theorem~\ref{thm:main} (The reducible case)}

\subsubsection{The reducible case: a unique maximal component}

Suppose $\X$ is reducible with $|I_{\max}|=1$, say $I_{\max}=\{1\}$, and
choose $J=[\alpha\beta]$ so that $\beta$ communicates with $C_1$ through every relevant chain in
$\mathcal P(i,j)$ containing no other maximal class.

\begin{proposition}
\label{prop:B}
Under this hypothesis, Theorem~\ref{thm:main} holds with $\X$'s Parry
measure replaced by the Parry measure $\mu_{C_1}$ of the unique maximal
recurrent component $\X_1$: for $t>t_c$, $P(t)=\xi$ and $\mu_{C_1}$ is the
unique equilibrium measure.
\end{proposition}

\begin{proof}

Since $I_{\max}=\{1\}$ is a singleton, every path in every relevant $\mathcal P(i,j)$ visits at most the single maximal class $C_1$, so
$d\equiv1$ along every relevant chain, exactly as in the mixing case.  By Lemma~\ref{lem:5.3},
$d_n\sim C\eta^n$, with $C$ depending on the entry/exit points into $C_1$. 
Equations~(6)--(8) will be same. Since $\varphi\equiv0$ on $\X$ and $h_\mu$ is upper semicontinuous, any
$\mu\in M_\sigma(\X)$ with $h_\mu=\xi$ is carried by $\bigcup_{i:\xi_i=\xi}
\X_i=\X_1$, and irreducibility of $C_1$ gives uniqueness of $\mu_{C_1}$ as
in Proposition~\ref{prop:A}.
\end{proof}

\begin{remark}
If instead $\beta,\alpha$ can only reach $C_1$ through a chain containing
another maximal class ($d\ge2$ for some, or every, $P\in\mathcal P(i,j)$), Lemma~\ref{lem:5.3} gives
$d_n\sim Cn^{d-1}\eta^n$, and $\lambda_{t,\xi}(\mathcal T_2)=\sum_{n\ge N}
n^{d-1-t}(Equation ~ (5))$ diverges for $t\le d$ rather than converging for all
$t\ge0$. This shifts the functional equation determining $t_c$ and is not
treated here; This is the reason we restrict our Main Theorem to non-chained case, where this issue does not arise.  
\end{remark}

\subsubsection{The reducible case: several maximal components}

Assume $\X$ has non-chained maximal classes (Definition~\ref{def:5.2})
with $|I_{\max}|=k\ge2$, say $I_{\max}=\{1,\dots,k\}$, and choose $J=[\alpha
\beta]$ so that $\beta,\alpha$ communicate with every $C_i$, $i\in I_{\max}$,
with no chain between any two distinct maximal classes, i.e, $d(\beta,\alpha)\le1$.

\begin{theorem}[General freezing transition]
\label{thm:general}
Let $\X$ be an SFT with entropy $\xi=\log\eta$ and non-chained maximal
classes $\X_1,\dots,\X_k$ ($k\ge1$), and $J$ as above. Let $\varphi$ be the
Hofbauer-type potential~(2). There is $t_c>0$ such that:
\begin{enumerate}
\item For $0<t<t_0$, there is a unique equilibrium measure $\mu_t$, of full
support; $P(t)$ is real analytic and $P(t)>\xi$.
\item For $t>t_0$, $P(t)=\xi$, and the set of equilibrium measures is
exactly the Parry simplex $\ParrySimplex$
\[
\Big\{\ \mu=\sum_{i=1}^k c_i\,\mu_{\X_i}\ :\ c_i\ge0,\ \textstyle\sum_i c_i
=1\ \Big\}.
\]
The transition is
freezing, and the equilibrium measure is unique (equal to the single Parry
measure of $\X$) if and only if $k=1$.
\end{enumerate}
\end{theorem}

\begin{proof}
By Lemma~\ref{lem:5.3} with $d\equiv1$ on every relevant chain,
each maximal class contributes $C_i\eta^n(1+o(1))$ to $d_n$, and by
disjointness of the families of return words routed through distinct
classes,
\[
d_n\sim\Big(\sum_{i=1}^k C_i\Big)\eta^n.
\]
Equations~(5)--(7) hold with $C=\sum_i C_i$, giving existence of $t_c$,
uniqueness and real-analyticity for $t<t_c$, and $P(t)=\xi$ for $t>t_c$
exactly as in the proof of Theorem~\ref{thm:main} in Section 4.

For the equilibrium measures at $t>t_c$: since $\varphi\equiv0$ on $\X$,
any $\mu\in M_\sigma(\X)$ with $h_\mu+t\int\varphi\,d\mu=\xi$ must (as in
Proposition~\ref{prop:B}, no equilibrium measure gives positive weight to
$J$ once $t>t_c$) have $h_\mu=\xi$ and support in $\X_{\max}=\bigcup_{i=1}^k
\X_i$. As the $\X_i$ are pairwise disjoint invariant closed subsystems,
$\mu$ decomposes as $\sum_i c_i\,\mu(\cdot\mid\X_i)$, $c_i=\mu(\X_i)$, and
by affinity of entropy,
\[
\xi=h_\mu=\sum_{i=1}^k c_i\, h_{\mu(\cdot\mid\X_i)}\le \sum_{i=1}^k c_i\,\xi_i
=\xi,
\]
forcing $h_{\mu(\cdot\mid\X_i)}=\xi_i=\xi$ whenever $c_i>0$. By intrinsic
ergodicity of each irreducible $\X_i$ \cite{Parry}, $\mu(\cdot\mid\X_i)=
\mu_{\X_i}$, so $\mu=\sum_i c_i\mu_{\X_i}$. Conversely every such convex
combination is invariant with entropy $\xi$ and $\int\varphi\,d\mu=0$,
hence is an equilibrium measure. This proves~(2); (1) is exactly
Propositions~\ref{prop:A}/\ref{prop:B} applied below the transition.
\end{proof}

\begin{remark}
Theorem~\ref{thm:general} specializes to Theorem~\ref{thm:main} when $\X$ is
mixing ($k=1$, trivially non-chained), and to Proposition~\ref{prop:A} when
$\X$ is irreducible but not mixing. The case $k\ge2$ exhibits a phenomenon
absent from the mixing setting: the freezing transition persists, but the
post-transition equilibrium state is no longer unique. Example~\ref{ex:5.5}
with two disjoint (non-chained) golden-mean shifts $A,B$ illustrates this
directly: for $t>t_c$, every measure $c\,\mu_A+(1-c)\mu_B$, $c\in[0,1]$, is
an equilibrium measure for $t\varphi$.
\end{remark}

\section{Conclusion}

We have shown that the Hofbauer-type potential defined via the distance to
a subshift of finite type $\X$ exhibits a freezing phase transition in full
generality. When $\X$ is mixing (Theorem~\ref{thm:main}) or, more generally,
irreducible (Proposition~\ref{prop:A}), the post-transition equilibrium
measure is unique and equal to the Parry measure of $\X$. When $\X$ is
reducible with several communicating classes tied for maximal entropy and
not chained to one another (Theorem~\ref{thm:general}), the transition
remains freezing, but the post-transition equilibrium measure ranges over
the full simplex spanned by the Parry measures of the maximal-entropy
components, a non-uniqueness phenomenon with no analogue in the
mixing or irreducible settings. The case of chained maximal-entropy
components, where the return-word counting function picks up a polynomial
correction $n^{d-1}$ (Lemma~\ref{lem:5.3}),
shifts the functional equation determining the transition point and is left
for future work.

\section*{AI Usage Disclaimer}
During the preparation of this manuscript, the author used an AI-based 
language assistant (Claude, Anthropic) to help refine the phrasing of the 
title and portions of the expository prose in the Introduction, and to 
suggest \LaTeX{} notation for the Parry simplex. All mathematical content, 
definitions, theorem statements, and proofs are entirely the author's own 
work. The author has reviewed, edited, and takes full responsibility for 
all content in this article.

\end{document}